\documentclass[11pt, amssymb,amsfonts]{amsart}

\usepackage{tikz}
\newcommand{\norm}[1]{ \left|  #1 \right| }
\newcommand{\Norm}[1]{ \left\|  #1 \right\| }

\def\P{{\hbox{\bf P}}}
\def\E{{\hbox{\bf E}}}

 at 10 true pt

\def\be#1{ \begin{equation}\label{#1} }

\def\bas{\begin{align*}}
\def\eas{\end{align*}}
\def\bi{\begin{itemize}}
\def\ei{\end{itemize}}

\def\emph#1{{\it #1}}
\def\textbf#1{{\bf #1}}

\theoremstyle{plain}
  \newtheorem{theorem}[subsubsection]{Theorem}

  \newtheorem{lemma}[subsubsection]{Lemma}
  \newtheorem{corollary}[subsubsection]{Corollary}

\theoremstyle{remark}
  \newtheorem{remark}[subsubsection]{Remark}

\theoremstyle{definition}

\include{psfig}

\title{Matrix perturbation bound via contour bootstrapping }
\author{Phuc Tran, Van Vu}

\thanks{Department of Mathematics, Yale University, 10 Hillhouse Ave, New Haven, Connecticut, USA.\\ \textit{Email:}\, \texttt{phuc.tran@yale.edu, van.vu@yale.edu}.}

\date{}
\begin{document}
\maketitle

\begin{abstract} Matrix perturbation bounds play an essential role in the design and analysis of spectral algorithms.
In this paper, we use a  "contour bootstrapping" argument to derive several new perturbation bounds. 
As applications, we discuss new bounds on the error occurring when one uses matrix sparsification to speed up 
the computation of spectral parameters. Another potential application is the estimation of the trade-off in computing with 
privacy. 

\end{abstract} 

\section{Introduction}
\subsection{Perturbation Theory} \label{section:classical} 

Let $A$ and $E$ be symmetric real matrices of size $n$, and set $\tilde A:= A+E$. We will view $A$ as the truth (or data) matrix, and $E$ as noise. 
Consider the spectral decomposition 
$$A =\sum_{i=1}^n \lambda_i u_i u_i^T , $$ where $\lambda_i$ are the eigenvalues of $u_i$ the corresponding eigenvectors.

We order the eigenvalues decreasingly,  $\lambda_1 \ge \lambda_2 \dots \geq \lambda_n$. We also order the singular values of $A$ in a similar manner 
$\sigma_1 \ge \sigma_2 \ge \dots \ge \sigma_n$. It is well known that
$\{\sigma_1, \dots, \sigma_n \} = \{| \lambda_1| , \dots,|  \lambda_ n| \} $. Thus we have $\sigma_i = | \lambda_{\pi (i) } | $ for some permutation 
$\pi \in S_n$. The gaps between the consecutive eigenvalues is defined as 
$$\delta_p := \lambda_p - \lambda_{p+1} . $$


 \begin{remark} To simplify the presentation, we assume that the eigenvalues (singular values) are different, so the eigenvectors (singular vectors) are well-defined (up to signs). However, 
 our results hold for matrices with multiple eigenvalues. 
 \end{remark}

Let $S$ be a subset of $\{1, \dots, n \}$ and $f=f(z)$ an analytic function.  
Consider the spectral decomposition 
$A= U \Lambda U^T $ where $\Lambda$ is the diagonal matrix, whose $i$th entry is $\lambda_i$, and $U$ is the matrix with the $i$th column being $u_i$. We define 

$$f_S(A)= U f(\Lambda) U^T ,$$  where $f(\Lambda)$ is the diagonal matrix whose $i$th entry is $f(\lambda_i)$ if $i \in S$ and zero otherwise. 

The functional $f_S (A)$ has been widely used in linear algebra; see \cite{Book1, HJBook}. We can represent many important 
operators with a proper choice of $f$ and $S$. For examples

\begin {itemize} 

\item For $f \equiv 1$, $f_S(A)=  \sum_{i \in S } u_i u_i^T $, the orthogonal projection onto the eigenspace spanned by the 
eigenvectors corresponding to the eigenvalues in with indices in $S$. An important subcase is the following:

\vskip2mm

\item Still with $f \equiv 1$, let $S= \{1, \dots, p \}$. Then  $f_S(A) $ is the orthogonal projection onto the eigenspace spanned by the leading 
$p$ eigenvectors. We will use the notation $\Pi_p$ for this operator. 

\vskip2mm

\item Let $S$ be a set such that the eigenvalues with indices in $S$ have the largest absolute values (in other words, the set 
$\{| \lambda_i |, i \in S \} $ is the set of the $p$ largest singular values of $A$. Then for $f(z)=z$, $f_S(A)$  is  the best rank $p$ approximation of $A$. We will use the notation $A_p$ for this operator. It is easy to see that if we arrange the eigenvalues decreasingly, then $S= \{1, \dots, k, n-p +(k+1), \dots, n\}$
for some $0 \le k \le p $. 

\end{itemize}

We are considering  perturbation bounds for   the difference
$$ \| f_S( A+ E) -f_S(A)  \| = \| f_S( \tilde A) - f_S(A) \| , $$
for an arbitrary pair of $f$ and $S$. 

In the next subsections, we will discuss two of the most important cases in applications: the eigenspace perturbation ($f_S(A) = \Pi_p$) and the perturbation of the best rank-$p$ approximation ($f_S(A) =A_p$). 

\subsection{Eigenspace Perturbation}
Consider a subset $S \subset \{1, \dots, n \}$ and denote by $\Pi _S$ the orthogonal projection onto the subspace spanned by the eigenvectors $u_i, i \in S$. 
Define $\delta_S =\min \{ | \lambda_i - \lambda_j | , i \in S, j \notin S \} $.  


\begin{theorem}[Davis-Kahan \cite{DKoriginal, Book1}] \label{DKgeneral} We have 
\begin{equation}    ||\tilde{\Pi}_S - \Pi_S || \leq \frac{2 \Norm{E}}{\delta_{S}}.
\end{equation}  

\end{theorem}

Consider the important special case where $S= \{1, \dots, p \}$,  we obtain the following estimate for the perturbation of $\Pi_p$, the orthogonal projection onto the leading $p$-eigenspace of $A$.

\begin{corollary}  \label{DKp} Let $\delta_p := \lambda_p -\lambda_{p+1}$. Then 
\begin{equation}
\Norm{\tilde{\Pi}_p - \Pi_{p}} \leq 2 \frac{||E||}{\delta_p}.
\end{equation}

\end{corollary}

\begin{remark}\label{gaptonoiseRM}  Notice that in any application of 
Theorem \ref{DKgeneral} or Corollary \ref{DKp}, one needs to assume 

\begin{equation} \label{gaptonoise} \delta_S \ge  \frac{2}{\epsilon} \| E \|, \end{equation} 

\noindent if we want to bound the perturbation by $\epsilon < 1$.  For instance, if we need the bound to be $o(1)$, then 
the gap $\delta_S$ must be larger than $\| E \| $ (the intensity of the noise) by an order of magnitude. It is well known that the Davis-Kahan bound is sharp up to a constant factor (see \cite{Book1, SS1, Kato1}). 
\end{remark}

\subsection{Low rank approximation}
Another frequently used result is the following estimate for the perturbation of low rank approximation. 

\begin {theorem}[Eckart-Young] \label{EY} 
\begin{equation}  \label{EKbound}  \| \tilde A_p - A _p \|  \le 2 ( \sigma_{p+1} + \| E \| ).  \end{equation} 

\end{theorem}

It is easy to derive this bound through the Eckart-Young theorem, which asserts that $ \| A - A_p \| = \sigma_{p+1} $. Indeed 
$$\| \tilde A_p - A _p \|  \le \| A - A_p\| + \| \tilde A - A \| + \| \tilde A - \tilde A_p \| \le \sigma_{p+1} + \|E \| + \tilde \sigma_{p+1} \le 2 (\sigma_{p+1} + \| E \| ), $$ where in the last estimate we use Weyl inequality. 
\begin{theorem}[Weyl \cite{We1, Book1}] For any $ 1 \leq i \leq n$, 
\begin{equation}
\norm{\tilde{\lambda}_i - \lambda_i} \leq ||E||.
\end{equation}
\end{theorem}

By this reason, we will refer to \eqref{EKbound} as  Eckart-Young bound.

 \begin{remark} If $ \max \left\lbrace ||E||, \sigma_p - \sigma_{p+1} \right\rbrace \leq \frac{\sigma_p}{2}$, then the RHS is $ 2(\sigma_p - (\sigma_p - \sigma_{p+1})+ ||E||)  \ge \|E \| + \sigma_p$.
\end{remark}
\begin{remark} If we want to bound  the perturbation in the  Frobenius norm, then we can use the fact that $\tilde A _p - A_p$ has rank at most $2p$ to deduce 
$$||\tilde{A}_p - A_p||_F \leq  \sqrt {2p} \| \tilde A_p - A_p \| .$$
\end{remark}

In this paper, we use the so-called ``contour bootstrapping"  argument to obtain new and improved  estimates, provided that there is a reasonable gap between the key eigenvalues of $A$. Our argument is general and can be applied to any function $f$, but as for applications, we will focus on the two settings discussed above (projection and low rank approximation). 

\vskip2mm

There have been numerous results which improve these classical 
perturbation bounds under different assumptions; see, for instance,  \cite{Book1, CCF1, CLX1, DOT1, EBW1, GTV1, HLMNV1, HJBook, JW1, KX1, MSZ1, PK1,  OVK 13, OVK 22, SS1, Sun1, Vu1, Wa1, Wei1, Wu1} and the references therein. 
In particular, the excellent books \cite{Book1, HJBook} contain many main developments up to 2013. 
To the best of our knowledge, our argument and results are new. We are going to illustrate this by focusing on 
the two special cases of projection and low rank approximation and comparing our results with the most recent ones in the literature.

The rest of the paper is organized as follows. In this next section, we describe our method, which leads to our key bound \eqref{keybound} in Lemma \ref{keylemma}. 

In Section \ref{newbounds}, we  use this bound to derive new results for the projection and low rank approximation 
problems. We present these new results concerning projection and low rank approximation in Subsections \ref{newDK} and \ref{newEY}, and then compare our new results with recent developments in the literature,  respectively.  In the rest of this section, we derive  general results for arbitrary function $f$. 

Randomly sparsifying the input matrix to speed up 
computation is a popular method in theoretical computer science. This is the core of RandNLA (randomized numerical linear algebra) see \cite{DM1} for a survey. Naturally, this speeding comes with a trade-off in the accuracy of the output, as the input matrix has changed. In Section \ref{applications}, we use our new bounds to estimate 
the error in several important algorithms.

The last section, Section \ref{details} is devoted to the technical details of the proofs of the results in Section \ref{newbounds}.


\section{The contour bootstrapping argument} \label{section:argument}

 We start with  Cauchy theorem. Consider a contour $\Gamma$ and a point $a \notin \Gamma$, together with a complex function $f$ which is analytic on $\Gamma$ and its interior. Then  
 \begin{equation} \label{Cauchy0}  \frac{1} {2 \pi {\bf i }} \int_{\Gamma}  \frac{f(z) }{z-a} dz = f(a) \end{equation} if  $a$ is inside $\Gamma$ and zero otherwise.  Now define 
 $$f (A)= U f(\Lambda ) U^T $$ where $f(\Lambda)$ is diagonal with entries $f(\lambda_i)$. Let  $\Gamma $ be a contour containing $\lambda_i, i \in S$ where $S$ is a subset of $\{1, \dots, n\}$, 
 and assume that all $\lambda_j, j \notin S$ are outside $\Gamma$. 
 We obtain the classical contour identity  
 \begin{equation} \label{contour-formula} 
  \frac{1} {2 \pi {\bf i }} \int_{\Gamma}   f (z) (z-A)^{ -1} dz  = \sum_{ i \in S} f(\lambda_i) u_i u_i ^T := f_S(A)  ; 
  \end{equation}

 \noindent  see \cite{Kato1, SS1}. By choosing an appropriate $f$, we  obtain many important matrix functionals

\begin{itemize} 

\item $f \equiv 1 $, $S= \{ 1, \dots, p \} $, then $\sum_{ i \in S} f(\lambda_i) u_i u_i ^T = \Pi_p$, the orthogonal projection on the subspace spanned by the first $k$ eigenvectors of $A$. 

\item$f \equiv 1$, $ S= \{ \pi (1), \dots, \pi (p) \}$, then $\sum_{ i \in S} f(\lambda_i) u_i u_i ^T = \Pi_{(p)}$, the orthogonal projection on the subspace spanned by the first $k$ singular vector of $A$.

\item $f(z) = z $, $ S= \{ \pi (1), \dots, \pi (p) \}$, then $\sum_{ i \in S} f(\lambda_i) u_i u_i ^T = \Pi_{(p)} A = A_p$, the best rank $p$ approximation of $A$. 

\item $f(z) = z^2  $, $ S= \{ \pi (1), \dots, \pi (p) \}$, then $\sum_{ i \in S} f(\lambda_i) u_i u_i ^T = \Pi_{(p)}  A^2 $, the best rank $p$ approximation of $A^2 $.

\end{itemize}  
 
For a moment, assume that the eigenvalues $\tilde \lambda_i ,  i\in S$ are inside $\Gamma$, and all $\tilde \lambda_j, j \notin S$ are outside. Then 

 \begin{equation} \label{contour-formula1} 
  \frac{1}{2 \pi {\bf i}}  \int_{\Gamma}  f (z) (z-\tilde A)^{-1} dz  = \sum_{ i \in S} f(\tilde \lambda_i) \tilde u_i  \tilde u_i ^T := f_S(\tilde A) .  \end{equation}  
\noindent This way, we derive a  contour identity for the perturbation 
 \begin{equation} \label{f_Sperturbation formula}
 f_S(\tilde A)- f_S(A)=  \frac{1} {2 \pi {\bf i} } \int_{\Gamma}  f (z) [(z-\tilde A)^{-1}- (z- A)^{-1} ]  dz. 
  \end{equation}   
 \noindent Now we bound the perturbation  by the corresponding integral 
  \begin{equation} \label{Boostrapinequality1}
  \Norm{f_S(\tilde A)- f_S(A)} \leq \frac{1}{2 \pi} \int_{\Gamma} \Norm{f (z) [(z-\tilde A)^{-1}- (z- A)^{-1} ]} dz.  
  \end{equation}
  
 \noindent  There have been attempts to control the RHS using series expansion and analytical tools; see for instance \cite[Part 2]{Kato1}. 
 We follow a different path, using a bootstrapping argument to reduce the estimation of the RHS of \eqref{Boostrapinequality1} to the estimation of a much simpler quantity, which can be computed directly.

  In what follows, we denote   $\frac{1}{2 \pi} \int_{\Gamma} \Norm{f (z) [(z-\tilde A)^{-1}- (z- A)^{-1} ]} dz$ by $F(f,S)$.
Using the resolvent formula
\begin{equation} M^{-1} - (M+N)^{-1}= (M+N)^{-1} N M^{-1} \end{equation} and the fact that $\tilde A =A+E$, we obtain 
\begin{equation} (z- A)^{-1} - (z-\tilde A)^{-1} =  (z-A)^{-1} E (z-\tilde A)^{-1} .  \end{equation}

\noindent Therefore, we can rewrite $F(f,S)$ as 
\begin{equation}
\begin{split}
F(f,S) & = \frac{1}{2 \pi} \int_{\Gamma} \Norm{ f(z) (z-A)^{-1} E (z- \tilde{A})^{-1}} dz \\
&= \frac{1}{2 \pi} \int_{\Gamma} \Norm{ f(z) (z-A)^{-1} E (z-A)^{-1}  -  f(z) (z-A)^{-1} E [(z-A)^{-1} - (z-\tilde{A})^{-1}]} dz.  
\end{split}
\end{equation}

\noindent Using triangle inequality, we obtain 
\begin{equation} \label{F(f,S)inequality1}
\begin{split}
F(f,S) & \leq \frac{1}{2 \pi} \int_{\Gamma} \Norm{f(z) (z-A)^{-1} E (z-A)^{-1}} dz + \frac{1}{2 \pi} \int_{\Gamma} \Norm{f(z) (z-A)^{-1} E [(z-A)^{-1} - (z-\tilde{A})^{-1}]}dz \\
&  \leq \frac{1}{2 \pi} \int_{\Gamma} \Norm{f(z) (z-A)^{-1} E (z-A)^{-1}} dz  + \\ &+ \frac{1}{2 \pi} \int_{\Gamma} \Norm{ (z-A)^{-1} E} \times \Norm{f(z) [(z-A)^{-1} - (z-\tilde{A})^{-1}]}dz \\
& \leq \frac{1}{2 \pi} \int_{\Gamma} \Norm{f(z) (z-A)^{-1} E (z-A)^{-1}} dz  + \\& +\frac{\max_{z \in \Gamma} \Norm{(z-A)^{-1} E} }{2 \pi} \int_{\Gamma} \Norm{f(z) [(z-A)^{-1} - (z-\tilde{A})^{-1}]}dz \\
& = \frac{1}{2 \pi} \int_{\Gamma} \Norm{f(z) (z-A)^{-1} E (z-A)^{-1}} dz  +\max_{z \in \Gamma} \Norm{(z-A)^{-1} E} \times F(f,S).
\end{split}
\end{equation}

Now we need our gap assumption. We assume that the distance between any eigenvalue in $S$ and any eigenvalue outside $S$ is at least $4  \| E \| $


\begin{equation} \label{keyassumption} \min _{ i\in S, j \notin S} |\lambda_i - \lambda_j |  \ge 4 \| E \| . \end{equation}

Under this assumption, we can draw  $\Gamma$ such that the distance from any eigenvalue of $A$ to $\Gamma$ is at least $2 \| E\| $. 
(The simplest way to do so is to construct $\Gamma$ out of horizontal and vertical segments, where the vertical ones bisect the intervals connecting an eigenvalue in $S$ with its nearest neighbor (left or right) outside $S$.) With such  $\Gamma$, we have 

$$ \max_{z \in \Gamma} \Norm{(z-A)^{-1} E} \leq  \frac{||E||}{2 \| E \| } = \frac{1}{2} \,(\text{by gap assumption}).$$

\noindent Together with \eqref{F(f,S)inequality1}, it follows that 
$$ F(f,S) \leq F_1(f,S) + \frac{1}{2} F(f,S),$$

\noindent where  
$$F_1(f,S):= \frac{1}{2 \pi} \int_{\Gamma} \Norm{f(z) (z-A)^{-1} E (z-A)^{-1}} dz.$$


Therefore, 
\begin{equation} \label{Boostrapinequality2}
 \frac{1}{2} F(f,S) \leq F_1(f,S).
\end{equation}

Notice that the gap assumption \eqref{keyassumption}  and Weyl's inequality ensure that $\tilde{\lambda}_i$ is inside the contour $\Gamma$ if and only if $i \in S$. 
Combining \eqref{Boostrapinequality1} and \eqref{Boostrapinequality2}, we obtain our key inequality

\begin{lemma} \label{keylemma} We have 
\begin{equation} \label{keybound}
\Norm{f_S(\tilde A)- f_S(A)} \leq 2 F_1(f,S).
\end{equation}

\end{lemma} 


\noindent This lemma is the heart of the bootstrapping argument. The remaining task is to compute/estimate $F_1(f, S)$. This task will vary from case to case, depending on $f$ and a proper choice for $\Gamma$. 
As illustration, in the next section, we demonstrate several new results derived by this method.

\section{New Perturbation Bounds } \label{newbounds}

\subsection{Eigenspace Perturbation}  \label{newDK}

Our first result is a refinement of Davis-Kahan theorem (Theorem \ref{DKgeneral}).  We fix an index $p$ and focus on the eigenspace spanned by the eigenvectors corresponding to 
the leading eigenvalues $\lambda_1, \dots, \lambda_p$. Notice that this space is well defined if $\lambda_p > \lambda_{p+1} $, even if there are indices  $i, j \le p$ where $\lambda_i = \lambda_j$.

\begin{theorem}[Eigenspace Perturbation] \label{main1}  Let $r \ge p$ be the smallest integer such that $ \frac{\norm{\lambda_p}}{2} \leq \norm{\lambda_p - \lambda_{r+1}} $, and set  $x: =\max_{i,j \le r} | u_i^T E u_j | $.
 Assume furthermore that  $$4||E|| \leq \delta_p: = \lambda_p - \lambda_{p+1}  \leq \frac{\norm{\lambda_p}}{4}. $$  Then 
\begin{equation}
\Norm{\tilde{\Pi}_p - \Pi_p} \leq 24 \left( \frac{||E||}{\norm{\lambda_p}} \log \left(\frac{6 \sigma_1}{\delta_p} \right) + \frac{r^2 x}{\delta_p} \right).
\end{equation}

\end{theorem}

To compare this theorem to Corollary \ref{DKp}, let us ignore the logarithmic term and $r^2$ (which is typically small in applications). Notice that in the first term on the RHS, we have $\| E \| $ in the numerator, but $\norm{\lambda_p} $ instead of $\delta _p =\lambda_p -\lambda_{p+1} $ in the denominator.  In the second term, the denominator is $\delta_p$, but the numerator is $x \le \| E \|$ (in particular, this can be much smaller than $\| E\| $ if $E$ is random, for instance). 
Indeed, let $E$ be a Wigner matrix (a random symmetric matrix whose upper diagonal entries are iid sub-gaussian random variables with mean 0 and variance 1). 

\begin{lemma}  Let $u$ and $v$ be two fixed unit vectors.  With probability $1-o(1)$, $\| E \| = (2+o(1)) \sqrt n $ and $u^T Ev = O(\log n)$. 
\end{lemma}

The first part of the lemma is a well-known fact in random matrix theory, while the second part is an easy consequence of Chernoff bound (the proof is left as an exercise; see \cite{ H1, St1, SS1}). One can easily extend this lemma to non-Wigner matrices. 
A crucial point here is that we can apply Theorem \ref{main1} and obtain a very good bound (say $o(1)$), even when the critical gap-to-noise condition \eqref{gaptonoise} (see the discussion in Remark \ref{gaptonoiseRM}), which is required in all previous applications of the original Davis-Kahan bound, fails.

\begin{corollary}  \label{cor:lowrank} Assume that $A$ has rank $r$ and $E$ is Wigner. Assume  furthermore. that $\delta_p \ge 2.01 \sqrt n$. Then with probability $1-o(1) $

$$ \Norm{\tilde{\Pi}_p - \Pi_p} = \tilde  O \left( \frac{||E|| }{\norm{\lambda_p}} + \frac{r^2 }{\delta_p} \right). $$ 

\end{corollary}

This is superior to Theorem \ref{DKp}, if $| \lambda_p | \gg \delta_p \gg r^2$. In the case when $E$ is gaussian, 
a similar result (under a weaker assumption) was obtained in \cite{OVK 13} by a different method. However, it seems highly non-trivial to extend the 
(fairly complicated) argument in this paper to general Wigner matrices.

{\bf \noindent Sharpness.}  We believe that under 
the given (gap) assumption, our bound is sharp, up to a logarithmic factor. Consider Theorem \ref{main1} and  assume that 
$r= O(1)$. In this case, we have 

$$ \| \tilde \Pi_p - \Pi_p \|  = \tilde O \left( \frac {\| E \| }{|\lambda_p| }  + \frac{x} {\delta_p} \right),$$  (where the notation 
$\tilde O $  hides the polylogarithmic term. We believe that both terms on the RHS are necessary. First, the noise-to-signal 
term $\frac{\| E \| }{ |\lambda_p |} $ has to be present, since if the intensity of the noise $\| E \| $  is much larger than 
the signal $|\lambda_p|$, then one does not expect to keep the eigenvectors from the data matrix. One can see a concrete example in \cite{B-GN1}, which is known as the BBP threshold phenomenon in random matrix theory. In particular, it is shown, under some assumption,  that if $E$ is a random gaussian matrix, and $\| E \| \ge |\lambda_1 | $, then $\tilde u_1$ is completely random. 
The second term $\frac{x}{ \delta_p }$ replaces the original bound $\frac{\| E \| }{\delta_p } $. The replacement of $\| E \|$  by $x$ is natural, as $x$ shows how $E$ interacts with the eigenvectors of $A$. We believe that this replacement is also optimal. 
In the random setting (Theorem \ref{main2}), $x$ is of order $\tilde O(1)$ and it seems unlikely that we can replace this with any smaller quantity. 

{\bf \noindent Comparison.} In most results concerning eigenspace perturbation that  we found in the literature,  both the assumption and the conclusion are different from ours, making the comparison impossible. The one that seems most relevant is Theorem 1 of \cite{JW1}, a very recent paper, which applies for the case when $A$ is positive semi-definite. Compared to this result, our advantage is that our assumption is very simple to check, the bound also has a simple form, and we do not require the restriction of positive semi-definiteness. Making a comparison regarding the strength of the two results is again very hard unless one restricts to a very specific setting. 

Finally, let us mention that there are many works for the case when  $E$ is random \cite{EBW1, HLMNV1, KX1, OVK 13, OVK 22, Vu1, Wu1, Z1}. In this setting, the most relevant works are, perhaps, \cite{OVK 22} and \cite{KX1}. 
In \cite{OVK 22}, the authors considered $E$ to be a random gaussian matrix (GOE), and basically proved the bound of Theorem \ref{main2}, even without the gap assumption. However, their argument uses special properties of gaussian matrix and is restricted to this setting. 
In  \cite{KX1}, the authors also considered $E$ being GOE.  Among others, they showed (under some assumption) that with probability $1-o(1),$
\begin{equation*}
\Norm{ \tilde{\Pi}_p - \Pi_p} = O \left[   \left(\frac{||E||}{\delta_{(p)}} \right)^2 + \frac{\sqrt{\log n}}{\delta_{(p)}} \right],
\end{equation*}
where $\delta_{(p)} = \min\{\delta_p, \delta_{p-1}\}.$

If we could ignore the term $\frac{\sqrt{\log n}}{\delta_{(p)}}$ on the RHS, then this is a quadratic improvement over the original Davis-Kahan theorem 
(Theorem \ref{DKp}).  Our bound is sharper than this bound if  $\frac{\delta_p^2}{\lambda_p } \leq ||E|| $. Again, it is not clear if the resutls in \cite{KX2} extends to general random matrices.

\subsection{Low rank approximation} \label{newEY} 

In this section, we derive two different bounds for the perturbation of best rank-$p$ approximation, each has an advantage in a certain setting. 
In this case, $f(z)=z$. 

Let $1 \leq k \leq p$ be a natural number such that 
$$\{ \lambda_{\pi(1)},...,\lambda_{\pi(p)}\}=\{ \lambda_1,...,\lambda_k > 0 \geq \lambda_{n-(p-k)+1},...,\lambda_n \}.$$
In other words, the $p$th singular value of $A$ is either $\lambda_k$ or  $| \lambda_{n- (p-k) +1}|$. We derive the ``halving distance"  $r$,  with respect to the two indices $p,k$,  as the smallest positive integer such that 

$$ \frac{\lambda_k}{2} \leq \lambda_k - \lambda_{r+1} \,\, {\rm  and } \,\, \frac{ |\lambda_{n-(p-k)+1}|} {2} \leq \lambda_{n-r+1 } - \lambda_{n-(p-k)+1}.$$

In the case when the singular values of $A$ decrease fast, $r$ is small. For instance, if $A$ is essentially low rank, then one can often set $r$ be the stable rank of $A$ (see the next section for the definition of stable rank) 
We set

$$\bar{x}:= \max_{ \substack{1 \leq i,j \leq r \\  n-r \le  i',j' \leq n }} \{  \norm{u_i^T E u_j}, \norm{ u_{i'}^T E u_{j'}} \}.$$

\begin{theorem}\label{main2.1} Assume that  $4||E|| \leq \delta_k \leq \frac{\lambda_k}{4}$ and $4||E|| \leq \delta_{n-(p-k)} \leq \frac{\norm{\lambda_{n-(p-k)+1}}}{4}$. Then 

\begin{equation}
\Norm{\tilde{A}_p -A_p} \leq 30 \left(||E||+r^2 \bar{x} \right) \left(\log \left(\frac{6 \sigma_1}{\delta_k} \right) + \log \left(\frac{6 \sigma_1}{\delta_{n-(p-k)}} \right) \right) + 30 r^2 \bar{x} \left( \frac{ \lambda_k}{\delta_k}+ \frac{ \norm{\lambda_{n-(p-k)+1}}}{\delta_{n-(p-k)}} \right).
\end{equation} 

\end{theorem}

In various applications, $r$ is small (order $\tilde O(1)$) and $x$ is negligible compared to $\| E \| $ (see the discussion following Theorem 3.2). In such settings, the  dominating term is $\| E \|$ and the bound becomes  $\tilde O (\| E \| )$. This improves 
the classical Eckart-Young bound (Theorem \ref{EKbound}) significantly in the case when the singular value $\sigma_{p+1} \gg \| E \| $.

\begin{corollary} \label{lowrankcor1} Assume that $A$ has rank $r$ and $E$ is a random symmetric matrix whose entries are independent random variables with mean zero and variance at most 1. Assume furthermore that with probability 1, all entries of 
$E$ are bounded by $K$. Then for any $1 \le p \le r$

$$\| \tilde A_p -A_p \| = \tilde O ( K \sqrt n ) . $$

\end{corollary}

Now we formulate our second theorem for low rank approximation.

\begin{theorem}\label{main2}  Define $1 \le k \le p$ as above. If $4||E|| \leq \min \{\delta_k, \delta_{n-(p-k)} \},$ then 
\begin{equation}
\Norm{\tilde{A}_p - A_p} \leq 6 ||E|| \left( \log \left( \frac{6 \sigma_1}{\delta_k} \right) + \frac{\lambda_k}{\delta_k} + \log \norm{\frac{6 \sigma_1}{\delta_{n-(p-k)}}} + \frac{\norm{\lambda_{n-(p-k)+1}}}{\delta_{n-(p-k)}}  \right).
\end{equation}
\end{theorem}

This theorem looks a bit simpler than the 
first, at the cost of having the term
$\frac{\lambda_k} {\delta_k}$  multiplied with $\| E\|$. (In the first theorem we only need to multiply $\log \frac{\lambda_k} {\delta_k} $. This result is superior 
to the Eckart-Young bound if the gap $\delta_p \gg \| E \| $.

\noindent {\bf Comparison.}  The only result we found having a similar (gap)  assumption to ours is Theorem 2.2 \cite{Vishnoi2}, with an application in the area of differential privacy. This theorem applies for $E$ being complex gaussian and its proof
 (similar to situations discussed in the previous subsection regarding \cite{OVK 22, KX1}) uses special properties of the complex gaussian matrix (GUE). One can directly derive a generalization of this theorem to any Wigner matrix using Theorem \ref{main2.1} and Theorem \ref{main2}. We will discuss more applications regarding computing with privacy in a forthcoming paper \cite{TVV}.



\subsection{ General $f$ } 

In this section, we consider a general function $f$. 
We fix $p$ and use short hand  $f_p(\tilde{A})$for  $f(\tilde{A}, \{1,2,...,p\})$. Similarly, we write 
$f_p (A)$ for  $ f(A, \{1,2,...,p\})$.

\begin{theorem}[Matrix perturbation] \label{main3}
Let $r$ be the smallest integer such that $\frac{\norm{\lambda_p}}{2} \leq \norm{\lambda_{r+1}- \lambda_p}$.  If $4||E|| \leq \delta_p \leq \frac{\norm{\lambda_p}}{4}$, then 
\begin{equation}
\Norm{f_p(\tilde{A}) - f_p(A)} \leq  24  \max_{z \in \Gamma_1} ||f(z)|| \times \left( \frac{||E|| }{\norm{\lambda_p}} \log \left(\frac{6 \sigma_1}{\delta_p} \right) + \frac{r^2 x}{\delta_p} \right) ,
\end{equation}
\noindent where  $\Gamma_1$ is the rectangle with vertices 
$$(x_0, T), (x_1, T), (x_1,-T), (x_0, -T),$$ 

\noindent with 
 $$x_0 := \lambda_p - \frac{\delta_p}{2},  x_1:= 2 \sigma_1, T:= 2\sigma_1,  \sigma_1 = \| A\|. $$

\end{theorem}

\section{Applications} \label{applications}

 \subsection{Stable rank of a matrix } \label{essrank} 
 
 Let $A$ be a matrix with singular values $\sigma_1 \ge \dots \ge \sigma_n$, we recall the definition of  the stable  rank of $A$ \cite[Section 7.6.1]{Vbook}  
 $$r_{stable} : = \frac{\sum_{i=1}^n \sigma_ i } {\sigma_1 } . $$ 
 
 It is clear that if $r \ge h r_{stable} $,  for some $h \ge 1$, then $\sigma_r \le h^{-1/2} \sigma_1  $.   The assumption that a data matrix has a small stable rank is often used in recent applications; see \cite[Chapter 7]{Vbook}, \cite{KL1, RV1, RV2}.

\subsection{Fast computation with random sparsification} 

An important and well-known method in computer science to speed up matrix calculation is to sparsify the 
input matrix (replacing a large part of the entries with zeroes).
  For instance, \cite{YZ1} shows that if a matrix is sparse, then matrix vector multiplication 
can be done much faster. Matrix-vector multiplication is the core operation of many algorithms. A good example is the standard iterative algorithm to compute the leading eigenvector \cite[Chapter 5]{NAbook}.  The survey \cite{DM1} discusses the (random) sparsification method in detail and contains a comprehensive list of references. 

One way to sparsify the input matrix is to zero out each entry with probability $1-\rho$, independently. Thus from the input matrix $A$, we get a sparse matrix $A'$ with density $\rho$. 
The matrix $\tilde A =\frac{1}{\rho} A'$ is a sparse random matrix, whose expectation is exactly $A$. Now assume that we are interested in computing $f_S(A)$. 
One then hopes that $f_S(\tilde A)$  would be a good approximation for $f_S(A)$. To summarize,  we gain on running time (thanks to the sparsification),  at the cost of 
the error term $\| f_S(\tilde A) -f_S(A) \| $.

Starting with the influential paper 
 \cite{AMc}, this procedure has been applied for a number of problems; 
 see, for example,  \cite{AFKMcS1, BDDMR1, BCN1, CLK1, CYZLK1, CmFmSq1, SKLZ1,SWZC1, SqCFT1, WMHZ1, XY1}.

We use our results to derive new estimates on the error term $\| f_S(\tilde A) -f_S(A) \|$. We provide two examples.

\subsubsection{Computing the leading eigenvector} 

Let $A$ be a symmetric matrix with bounded entries $\| A \|_{\infty} \le K$. We want to estimate the first eigenvector of $A$.  
We sparsify $A$ as above, keeping the symmetry (i.e., zeroing out a pair of entries at a time). 
The above setting gives $\tilde A = A +E$, where $E$ is a random symmetric with entries bounded by $K/\rho$ and zero mean. We use  the following estimates for random matrices: 
\begin{itemize}
\item  By \cite[Theorem 1.4 ]{Vu1}, if $\rho \gg \frac{\log^4 n}{n}$, then $||E|| \leq 2K \sqrt{n/\rho}$ almost surely.
\item By \cite[Lemma 35]{OVK 13}, for any fixed unit vectors $u,v$, if $K/\rho > 1$, then 
$$\P \left(|u^T E v| \geq  \frac{ 2\sqrt{2} K}{\sqrt{\rho}} t \right) \leq \exp \left( -t^2 \right).$$

\noindent This implies, by the union bound, that 
$$\P \left( x \geq  \frac{2 \sqrt{2} K}{\sqrt{\rho}} \log n \right) \leq r^2 n^{-2}, $$
where we define  $x: =\max_{i,j \le r} | u_i^T E u_j | $.
\end{itemize}

\noindent Using these technical facts, it is easy to see that  Theorem \ref{main1} yields

\begin{corollary}  \label{sparse1} 
Let $r \ge 1$ be the smallest integer such that $ \frac{\norm{\lambda_1}}{2} \leq \norm{\lambda_1 - \lambda_{r+1}} $. Let $\rho$ be a positive number such that $1 \geq \rho \gg \frac{\log^4 n}{n}$. Assume furthermore that  $$8K \sqrt{n/\rho} \leq \delta_1 := \lambda_1 - \lambda_{2}  \leq \frac{\norm{\lambda_1}}{4}. $$  Then with high probability, 
\begin{equation} \label{sparse11} 
\Norm{\tilde{u}_1 - u_1} \leq \frac{72K}{\sqrt{\rho}} \left( \frac{ \sqrt{n}}{\norm{\lambda_1}} \log \left(\frac{6 \sigma_1}{\delta_1} \right) + \frac{ r^2 \log n }{ \delta_1} \right).
\end{equation}

\end{corollary} 

By the remark in subsection \ref{essrank}, it is clear that we can take $r= 4 r_{stable} (A) $. Corollary \ref{sparse1} is efficient in the case the stable rank is small. 
To make a comparison, let us notice that if we use Davis-Kahan bound (combined with the above estimate on $\| E \| $), we would obtain 

\begin{equation} \label{sparse12}  \Norm{\tilde{u}_1 - u_1} = O \left( \frac{\| E \| }{\delta_1} \right) = O \left( \frac{K \sqrt {n/ \rho} } {\delta_1 } \right).  \end{equation} 

To compare with Corollary \ref{sparse1}, notice that the first term on the RHS of \eqref{sparse11} is superior to the RHS of \eqref{sparse12} if $\lambda_1 \gg \delta_1$, which occurs if the two leading eigenvalues are relatively 
close to each other. The second term in the RHS of \eqref{sparse11} is also small compared to the RHS of \eqref{sparse12} if the stable rank is small, say $r_{stable} (A)= o(  \sqrt {  n^{1/2}  \log ^{-1} n } )= o( n^{1/4 } \log^{-1/2} n )$.

\subsubsection{Computing low rank approximation}

The main goal of the influential paper \cite{AMc} is to compute the low rank approximation $A_p$, for a  small $p$. 
In this paper, the authors bounded the difference between $\Norm{A -\tilde{A}_p}$ and $\Norm{A - A_p}$.  They proved

\begin{theorem} \cite{AMc}  If $\rho \gg \frac{\log^6 n}{n}$, then with probability at least $ 1 -\frac{1}{n}$, 
\begin{equation} \label{AMcbound}
\Norm{A -\tilde{A}_p} \leq \Norm{A - A_p} + 14 K \sqrt{n/\rho}.
\end{equation}

\end{theorem}

On the other hand, by using Theorem \ref{main2.1} and Theorem \ref{main2}, we can directly bound the difference $\Norm{\tilde{A}_p -A_p}$. This enables one to gain information about $A_p$, knowing 
$\tilde A_p$, which would be important in practice.   From Theorem \ref{main2} , we  derive  the following corollary






\begin{corollary} \label{lowrankcor} 
Define $1 \le k \le p$ and $\rho$ as above. If $ 8K \sqrt{n/\rho} \leq \min \{\delta_k, \delta_{n-(p-k)} \},$ then with probability at least $1 -\frac{r^2}{n^2}$,
\begin{equation}
\Norm{\tilde{A}_p - A_p} \leq 12 K \sqrt{n/\rho} \left( \log \left( \frac{6 \sigma_1}{\delta_k} \right) + \frac{\lambda_k}{\delta_k} + \log \norm{\frac{6 \sigma_1}{\delta_{n-(p-k)}}} + \frac{\norm{\lambda_{n-(p-k)+1}}}{\delta_{n-(p-k)}}  \right).
\end{equation}

\end{corollary}


To analyze  this new result, it is simpler to 
 assume that $A$ is positive definite. In this case,  the eigenvalues are the singular values, and  the  RHS simplifies to 
\begin{equation} \label{posdef2cor}  12 K \sqrt{n/\rho} \left(\log \left( \frac{6 \sigma_1}{\delta_p} \right)+ \frac{ \sigma_p}{\delta_p} \right) =\tilde O( K\sqrt {n\rho^{-1} } \frac{\sigma_p} {\delta_p}).  \end{equation}

\noindent Let us compare this with the Eckart-Young bound (Theorem \ref{EY}), which gives the bound 

$$O (\| E \|  + \sigma_{ p+1}  ) = O( K \sqrt { n \rho^{-1}  } + \sigma_{p+1}) . $$

The main difference here is that we have  $ K \sqrt { n \rho^{-1} }  \frac{ \sigma_p}{\delta_p}$ instead of  $\sigma _{p+1} $. This works in our favor if the gap $\delta_p$ is sufficiently large compared to 
$K \sqrt {n \rho^{-1} }$. First, notice that this is already part of our assumption of Corollary \ref{lowrankcor}.  More importantly, this happens often in practice, as to achieve a good low rank approximation, one naturally chooses a place where the gap is large. 
Quantitatively,  if $\sigma_{p+1} = c \sigma_p = \omega (K \sqrt {n \rho^{-1} } )$, for some constant $c<1$, then 
$\delta_p = (1-c) \sigma_p =\omega (K \sqrt {n/\rho} )$, and our bound is superior.  One can make the same comparison for the general case when $A$ is not positive definite, but the details are a bit more involved.  
 


\section{Proofs of the main theorems} \label{section: proof}

All the proofs are essentially bounding the quantity $F_1(f,S) $ discussed in Section \ref{section:argument}.

\subsection{Proof of Theorem \ref{main1} }
Recall that 
$$2 \pi F_1(1,\{1,2,...,p \}):= \int_{\Gamma} \Norm{(z-A)^{-1} E (z-A)^{-1}} dz,$$
here the contour $\Gamma$ is a rectangle, whose  four vertices are 
$$(x_0, T), (x_1, T), (x_1,-T), (x_0, -T),$$
with 
 $$x_0 := \lambda_p - \delta_p/2, x_1:= 2 \lambda_1, T:= 2\lambda_1.$$

Now, we split $\Gamma$ into four segments: 
\begin{itemize}
\item $\Gamma_1:= \{ (x_0,t) | -T \leq t \leq T\}$. 
\item $\Gamma_2:= \{ (x, T) | x_0 \leq x \leq x_1\}$. 
\item $\Gamma_3:= \{ (x_1, t) | T \geq t \geq -T\}$. 
\item $\Gamma_4:= \{ (x, -T)| x_1 \geq x \geq x_0\}$. 
\end{itemize}
\usetikzlibrary{decorations.pathreplacing}
$$\begin{tikzpicture}

\coordinate (A) at (4,0);
\node[below] at (A){$\lambda_{p+1}$};
\coordinate (A') at (6, 0);
\node[below] at (A'){$\lambda_p$};
\coordinate (B) at (5,0);
\node[above left] at (B){$\Gamma_1$};
\coordinate (C) at (8,0);
\node[below] at (C){$\lambda_1$};
\coordinate (D) at (125mm,0);
\node[above right] at (D){$\Gamma_3$};
\coordinate (E) at (13,0);
\coordinate(B') at (5,5mm);
\coordinate (F) at (9,1);
\node[above] at (F){$\Gamma_2$};
\coordinate (G) at (9,-1);
\node[below] at (G){$\Gamma_4$};
\draw (0,0) -- (A);
\draw[very thick,blue] (A) -- (A');
\draw[very thick, brown] (A')  -- (C);
\draw (C) -- (D) -- (E);
\draw[->] (B) -- (B');
\draw (B') -- (5,1);

\draw [->] (5,1) -- (9,1);
\draw (9,1) -- (125mm,1) -- (D);
\draw [->] (D) -- (125mm, -5mm);
\draw (125mm,-5mm) -- (125mm,-1);
\draw [->] (125mm,-1) -- (9,-1); 
\draw (9,-1) -- (5,-1) -- (B);

\end{tikzpicture}.$$

Therefore, 
$$2\pi F_1 = \sum_{k} M_k, M_k:= \int_{\Gamma_k} \Norm{(z-A)^{-1} E (z-A)^{-1} } dz.$$
We will use the following lemmas to bound $M_1, M_2, M_3, M_4$ from above. The proofs of these lemmas are delayed to the next section.

\begin{lemma} \label{F1f1M1bound}
Under the assumption of Theorem \ref{main1}, 
\begin{equation}
M_1 \leq 70 \left( \frac{||E|| }{\norm{\lambda_p}} \log \left(\frac{6 \sigma_1}{\delta_p} \right) + \frac{r^2 x}{\delta_p} \right) .
\end{equation}
\end{lemma}
\begin{lemma} \label{F1f1M2bound}
Under the assumption of Theorem \ref{main1}, 
\begin{equation}
M_2, M_4 \leq \frac{||E|| \cdot \norm{x_1-x_0}}{T^2}. 
\end{equation}
\end{lemma}
\begin{lemma} \label{F1f1M3bound}
Under the assumption of Theorem \ref{main1}, 
\begin{equation}
M_3 \leq \frac{4||E||}{\norm{x_1-\lambda_1}}.
\end{equation}
\end{lemma}
By our setting of $x_1,x_0,T$,  we have 
$$M_2,M_4 \leq \frac{||E|| 2\sigma_1}{4 \sigma_1^2}= \frac{||E||}{2 \sigma_1}, M_3 \leq \frac{4||E||}{\sigma_1}.$$
Therefore, 
\begin{equation} 
\begin{split}
F_1 & \leq \frac{1}{2\pi}\left( M_1+M_2+M_3+M_4 \right) \\
& \leq  \frac{70}{2\pi} \left( \frac{||E||}{\norm{\lambda_p}} \log \left(\frac{6 \sigma_1}{\delta_p} \right) + \frac{r^2 x}{\delta_p} \right)  + \frac{||E||}{ 2\pi \sigma_1} + \frac{4||E||}{2 \pi \sigma_1} \\
& \leq \frac{72}{2\pi} \left( \frac{||E||}{\norm{\lambda_p}} \log \left(\frac{6 \sigma_1}{\delta_p} \right) + \frac{r^2 x}{\delta_p} \right) \,\,\,(\text{since}\,\, \log \left(\frac{6 \sigma_1}{\delta_p} \right) > \log 24 > 3) \\
& \leq 12 \left( \frac{||E||}{\norm{\lambda_p}} \log \left(\frac{6 \sigma_1}{\delta_p} \right) + \frac{r^2 x}{\delta_p} \right). 
\end{split}
\end{equation}
Since $\Norm{\tilde{\Pi}_p -\Pi_p} \leq 2 F_1$, we finally obtain 

$$\Norm{\tilde{\Pi}_p -\Pi_p}  \leq 24 \left( \frac{||E||}{\norm{\lambda_p}} \log \left(\frac{6 \sigma_1}{\delta_p} \right) + \frac{r^2 x}{\delta_p} \right).$$

\subsection{Proof of Theorem \ref{main3}} Notice that 
\begin{equation}
\begin{split}
F_1(f,\{1,2,...,p\}) \leq \max_{z \in \Gamma}||f(z)|| \times F_1(1,\{1,2,...,p\}).
\end{split}
\end{equation}
Therefore, applying the result from previous subsection, we conclude that with high probability
\begin{equation}
\begin{split}
\Norm{f_p(\tilde{A}) - f_p(A)} \leq F(f,\{1,2,...,p\}) & \leq 2 F_1(f,\{1,2,...,p\})\\
& \leq 2 \max_{z \in \Gamma}||f(z)|| \times F_1(1,\{1,2,...,p\}) \\
& \leq 24 \max_{z \in \Gamma}||f(z)|| \times  \left( \frac{||E||}{\norm{\lambda_p}} \log \left(\frac{6 \sigma_1}{\delta_p} \right) + \frac{r^2 x}{\delta_p} \right).
\end{split}
\end{equation}
We complete our proof. 

\subsection{Proof of Theorem \ref{main2} }

 Let $1 \leq k \leq p$ be a natural number such that 
$$\{ \lambda_{\pi(1)}, \lambda_{\pi(2)}, ..., \lambda_{\pi(p)} \} = \{ \lambda_1, \lambda_2,...,\lambda_k > 0 \geq \lambda_{n-(p-k)+1}, \lambda_{n-(p-k)+2},...,\lambda_{n} \}.$$ 
Thus, we can split $A_p$ as  $ A_k + B_{p-k},$ in which 
$$B_{p-k} = \sum_{n \geq i \geq n-(p-k)+1} \lambda_{i} u_i u_i^T.$$
Similarly, $\tilde{A}_p = \tilde{A}_k + \tilde{B}_{p-k}$. Therefore, 
$$ \Norm{\tilde{A}_p - A_p} = \Norm{\tilde{A}_k + \tilde{B}_{p-k} -A_k - B_{p-k}} \leq \Norm{\tilde{A}_k - A_k} + \Norm{\tilde{B}_{p-k} - B_{p-k}}.$$
Applying the contour bootstrapping argument on $\Norm{\tilde{A}_k - A_k}$ with contour $\Gamma^{[1]}$ and on $ \Norm{\tilde{B}_{p-k} - B_{p-k}}$ with another contour $\Gamma^{[2]}$ (we define these contours later), we obtain 
\begin{equation} \label{splitF_1f2}
\begin{split}
& \Norm{\tilde{A}_k - A_k} \leq 2 F_1^{[1]}:=  \frac{1}{2\pi} \int_{\Gamma^{[1]}} \Norm{z (z-A)^{-1} E (z-A)^{-1}} dz, \\
& \Norm{\tilde{B}_{p-k} - B_{p-k}} \leq 2 F_1^{[2]}:= \frac{1}{2\pi} \int_{\Gamma^{[2]}} \Norm{ z (z-A)^{-1} E (z-A)^{-1}} dz.
\end{split}
\end{equation}
Here, $\Gamma^{[1]}$ and $ \Gamma^{[2]}$ are rectangles, whose vertices are 
$$\Gamma^{[1]}: (a_0, T), (a_1, T), (a_1,-T), (a_0, -T),$$
with 
 $$a_0 := \lambda_k - \delta_k/2, a_1:= 2 \lambda_1, T:= 2 \sigma_1;$$
and 
$$\Gamma^{[2]}: (b_0, T), (b_1, T), (b_1,-T), (b_0, -T),$$
with 
 $$b_0 := \lambda_{n-(p-k)+1} + \delta_{n-(p-k)}/2, b_1:= 2 \lambda_n, T:= 2 \sigma_1.$$

Now, we are going to bound $F_1^{[1]}$. First, we split $\Gamma^{[1]}$ into four segments: 
\begin{itemize}
\item $\Gamma_1:= \{ (a_0,t) | -T \leq t \leq T\}$. 
\item $\Gamma_2:= \{ (x, T) | a_0 \leq x \leq a_1\}$. 
\item $\Gamma_3:= \{ (a_1, t) | T \geq t \geq -T\}$. 
\item $\Gamma_4:= \{ (x, -T)| a_1 \geq x \geq a_0\}$. 
\end{itemize} 
Therefore, 
$$F_1^{[1]} = \sum_{l=1}^{4} \int_{\Gamma_l}  \Norm{ z (z-A)^{-1} E (z-A)^{-1}} dz.$$
Notice that 
$$\Norm{z (z-A)^{-1} E (z-A)^{-1}} \leq ||E|| \frac{|z|}{\min_i |z-\lambda_i|^2},$$
we further obtain
\begin{equation} \label{F_1f2inequality1}
 2 \pi F_1^{[1]} \leq ||E|| \left(\sum_{l=1}^4 N_l \right),
\end{equation}
in which 

$$N_l:= \int_{\Gamma_l}  \frac{|z|}{\min_i |z-\lambda_i|^2} dz \,\,\, \text{for}\,\, l=\overline{1,4}. $$ 
We use the following lemmas (again the  proofs are delayed to the next section). 

\begin{lemma} \label{lemma: N_1,3} Under our contour assumption, 

\begin{equation}
\begin{split}
& N_1 \leq \frac{8 a_0}{\delta_k} + 4 \log \norm{\frac{3T}{\delta_k}},\\
& N_3 \leq \frac{4 a_1}{ (a_1-\lambda_1)} + 4  \log \norm{\frac{3T}{(a_1 -\lambda_1)}}.
\end{split}
\end{equation}
\end{lemma}

\begin{lemma} \label{lemma: N2,4} Under our contour assumption, 
\begin{equation}
N_2, N_4 \leq \frac{\sqrt{2}(a_1-a_0)}{T}.
\end{equation}
\end{lemma}
By our setting of $a_0, a_1, T$, we have 
$$N_1 \leq \frac{8 \lambda_k}{\delta_k} + 4\log \norm{\frac{6 \sigma_1}{\delta_k}} , N_3 \leq 8 + 4 \log 6,\,\,\text{and}\,\, N_2, N_4 \leq \sqrt{2}.$$
Therefore, 
\begin{equation} \label{F_1'bound}
\begin{split}
F_1^{[1]}&  \leq  \frac{||E||}{2 \pi } \left(8+ 4 \log 6+ 2\sqrt{2} + \frac{8 \lambda_k}{\delta_k} + 4  \log \norm{\frac{6 \sigma_1}{\delta_k}} \right) \\
& \leq  \frac{||E||}{2\pi} \left( 15 \log \norm{\frac{6 \sigma_1}{\delta_k}} + \frac{8 \lambda_k}{\delta_k} \right) \\
& \leq 3 ||E|| \left(  \log \norm{\frac{6 \sigma_1}{\delta_k}} + \frac{ \lambda_k}{\delta_k} \right)
\end{split}
\end{equation}

Applying a similar argument on contour $\Gamma^{[2]}$, we obtain 
\begin{equation} \label{F_1"bound}
F_1^{[2]} \leq 3 ||E|| \left( \log \norm{\frac{6 \sigma_1}{\delta_{n-(p-k)}}} + \frac{\norm{\lambda_{n-(p-k)+1}}}{\delta_{n-(p-k)}} \right).
\end{equation}

Combining \eqref{splitF_1f2}, \eqref{F_1'bound} and \eqref{F_1"bound}, we complete our proof. 
\subsection{Proof of Theorem \ref{main2.1}} Intuitively, in this subsection, we will combine the arguments from the proofs of Theorems \ref{main1} and \ref{main2}. First, we still split $(\tilde{A}_p,A_p)$ into $(A_k,B_{p-k}, \tilde{A}_k, \tilde{B}_{p-k}$ and apply the contour bootstrapping argument on $\Norm{\tilde{A}_k -A_k}, \Norm{\tilde{B}_{p-k} -B_{p-k}}$. We obtain 
$$\Norm{\tilde{A}_p -A_p} \leq  2\left( F_1^{[1]} + F_1^{[2]} \right).$$
However, we will treat $F_1^{[1]}, F_1^{[2]}$  a bit differently. Indeed, 
$$2\pi F_1^{[1]} \leq M_1 + ||E|| \left(N_2+N_3 +N_4 \right),$$
in which 
$$M_1:= \int_{\Gamma_1} \Norm{ \sum_{i,j} \frac{z}{(z-\lambda_i)(z-\lambda_j)} u_i u_i E u_j u_j^T} dz.$$

We use the following lemma. 
\begin{lemma} \label{lemma: M_1f2} Under the assumption of Theorem \ref{main2.1}, 
\begin{equation}
M_1 \leq r^2 \bar{x} \left( \frac{ 8 \lambda_k}{\delta_k} + 2\log \left(\frac{3T}{\delta_k} \right)   \right) + 80 ||E||\log \left(\frac{3T}{\delta_k}  \right).
\end{equation}
\end{lemma}

Using above lemma and the upper bounds on $N_2,N_3,N_4$ from the previous subsection, we obtain 
\begin{equation}
\begin{split}
2 \pi F_1^{[1]} & \leq  r^2 \bar{x} \left( \frac{8 \lambda_k}{\delta_k}+ 2 \log \left(\frac{3T}{\delta_k}  \right) \right) + 80 ||E|| \log \left(\frac{3T}{\delta_k} \right) + ||E|| \left(2\sqrt{2}+ 8+ 4\log 6 \right) \\
& \leq  r^2 \bar{x} \left( \frac{8 \lambda_k}{\delta_k}+ 2 \log \left(\frac{6 \sigma_1}{\delta_k}  \right) \right) + 80 ||E|| \log \left(\frac{6 \sigma_1}{\delta_k} \right) + 11 ||E|| \log \left(\frac{6 \sigma_1}{\delta_k} \right)   \\
& = 91 \left(||E|| \log \left(\frac{6 \sigma_1}{\delta_k} \right) + r^2 \bar{x} \frac{ \lambda_k}{\delta_k} + r^2 \bar x  \log \left(\frac{6 \sigma_1}{\delta_k}  \right) \right).
\end{split}
\end{equation}
Thus, 
$$F_1^{[1]} \leq 15 \left(||E|| \log \left(\frac{6 \sigma_1}{\delta_k} \right) + r^2 \bar{x} \frac{ \lambda_k}{\delta_k} + r^2 \bar x  \log \left(\frac{6 \sigma_1}{\delta_k}  \right) \right).$$
Similarly, 
$$F_1^{[2]} \leq 15 \left(||E|| \log \left(\frac{6 \sigma_1}{\delta_{n-(p-k)}} \right) + r^2 \bar{x} \frac{ \norm{\lambda_{n-(p-k)+1}}}{\delta_{n-(p-k)}} + r^2 \bar x  \log \left(\frac{6 \sigma_1}{\delta_{n-(p-k)}}  \right) \right).$$
Therefore, 
$$\Norm{\tilde{A}_p -A_p} \leq 30  \left(||E|| \log \left( \frac{36 \sigma_1^2}{\delta_k \delta_{n-(p-k)}} \right)  + r^2 \bar x \frac{ \lambda_k}{\delta_k} + r^2 \bar x \frac{ \norm{\lambda_{n-(p-k)+1}}}{\delta_{n-(p-k)}}+ r^2 \bar x \log \left( \frac{36 \sigma_1^2}{\delta_k \delta_{n-(p-k)}} \right)  \right).$$

\section{Proofs of the lemmas } \label{details}

\subsection{Proofs of Lemma \ref{F1f1M2bound} and Lemma \ref{F1f1M3bound}}
Notice that 
$$\Norm{(z-A)^{-1} E (z-A)^{-1}} \leq \frac{||E||}{\min_{i}|z- \lambda_i|^2}.$$
Therefore, 
\begin{equation} \label{F_1f1inequality1}
M_2 \leq \int_{\Gamma_2}  \frac{1}{\min |z-\lambda_i|^2 } ||E|| dz = \Norm{E} \int_{\Gamma_2}  \frac{1}{\min |z-\lambda_i|^2 } dz.
\end{equation}
Moreover, 
\begin{equation} \label{F_1f1inequality3}
\begin{split}
& \int_{\Gamma_2}   \frac{1}{\min |z-\lambda_i|^2 } dz \\
& = \int_{\Gamma_2}  \frac{1}{ \min_i ((x-\lambda_i)^2+T^2)} dx \\
& \leq \int_{\Gamma_2} \frac{1}{T^2} dx \\
& = \frac{|x_1 - x_0|}{T^2}.
\end{split}
\end{equation}
Therefore, $M_2 \leq \frac{||E|| \cdot \norm{x_1-x_0}}{T^2}.$
Similarly, we also obtain that $M_4 \leq \frac{||E|| \cdot \norm{x_1-x_0}}{T^2}.$\\
Next, 
\begin{equation} \label{F_1f1inequality4}
\begin{split}
M_3 & \leq  ||E|| \int_{\Gamma_3}   \frac{1}{\min |z-\lambda_i|^2 } dz \\
& = ||E|| \int_{\Gamma_3}  \frac{1}{ \min ((x_1 - \lambda_i)^2 + t^2)} dt \\
& \leq ||E|| \int_{-T}^{T} \frac{1}{t^2 + (x_1-\lambda_1)^2} dt \\
& \leq \frac{ 4 ||E||}{x_1-\lambda_1}.
\end{split}
\end{equation}
The last inequality is true by the following lemma. 
\begin{lemma} \label{ingegralcomputation1}
Let $a, T$ be positive numbers such that $a \leq T$. Then, 
\begin{equation}
\int_{-T}^{T} \frac{1}{t^2+a^2} dt \leq \frac{4}{a}.
\end{equation}
\end{lemma}
\begin{proof}
We have 
\begin{equation}
\begin{split}
\int_{-T}^{T} \frac{1}{t^2+a^2} dt & = \int_{-T}^{T} \frac{1}{a^2 u^2 +a^2} dt (\,t = au) \\
& = \frac{1}{a} \int_{-T/a}^{T/a} \frac{1}{u^2+1} du \\
& = \frac{2}{a} \int_{0}^{T/a} \frac{1}{u^2+1} du \\
& \leq \frac{2}{a} \left( \int_{0}^{1} \frac{1}{u^2+1} du + \int_{1}^{T/a} \frac{1}{u^2} du \right) \\
& \leq \frac{2}{a} \left( 1 + 1 - \frac{a}{T} \right) \leq \frac{4}{a}.
\end{split}
\end{equation}
\end{proof}

\subsection{Proof of Lemma \ref{F1f1M1bound}}
Using the spectral  decomposition  $(z-A)^{-1} = \sum_{i} \frac{u_i u_i^T}{(z- \lambda_i)}$, we can rewrite $M_1$ as 
\begin{equation}
M_1 = \int_{\Gamma_1} \Norm{ \sum_{i,j} \frac{1}{(z-\lambda_i)(z-\lambda_j)} u_i u_i^T E u_j u_j^T} dz.
\end{equation}

\noindent Recall that $x:=\max_{1 \leq i,j \leq r} \norm{u_i^T E u_j}$. Using the  triangle inequality, we have 
\begin{equation}
\begin{split}
M_1 & \leq \int_{\Gamma_1} \Norm{ \sum_{i,j \leq r} \frac{1}{(z-\lambda_i)(z-\lambda_j)} u_i u_i^T E u_j u_j^T} dz + \int_{\Gamma_1} \Norm{ \sum_{i,j > r} \frac{1}{(z-\lambda_i)(z-\lambda_j)} u_i u_i^T E u_j u_j^T} dz \\
& + \int_{\Gamma_1} \Norm{ \sum_{\substack{i \leq r < j \\ i > r \geq j}} \frac{1}{(z-\lambda_i)(z-\lambda_j)} u_i u_i^T E u_j u_j^T} dz.
\end{split}
\end{equation}

\noindent Consider the first term 
\begin{equation}
\begin{split}
\int_{\Gamma_1} \Norm{ \sum_{i,j \leq r} \frac{1}{(z-\lambda_i)(z-\lambda_j)} u_i u_i^T E u_j u_j^T} dz & \leq \sum_{i,j \leq r} x \int_{\Gamma_1} \frac{1}{\sqrt{((x_0 -\lambda_i)^2+t^2)((x_0 -\lambda_j)^2+t^2)}} dt \\
& \leq r^2 x \int_{-T}^{T} \frac{1}{t^2 +(\delta_p/2)^2} dt \,(\text{since}\,\,|x_0 - \lambda_i| \geq \delta_p/2 \,\forall i) \\
& \leq \frac{8 r^2 x}{\delta_p} \,\,\, (\text{by Lemma \ref{ingegralcomputation1}}).
\end{split}
\end{equation}

Next, we apply the argument for bounding $M_2$  to bound  the second term,

\begin{equation}
\begin{split}
\int_{\Gamma_1} \Norm{ \sum_{i,j > r} \frac{1}{(z-\lambda_i)(z-\lambda_j)} u_i u_i^T E u_j u_j^T} dz & = \int_{\Gamma_1} \Norm{ \left(\sum_{i >r} \frac{u_iu_i^T}{z- \lambda_i} \right) E \left( \sum_{i > r} \frac{u_i u_i^T}{z -\lambda_i} \right)} dz \\
& \leq \int_{\Gamma_1} \Norm{\left(\sum_{i >r} \frac{u_iu_i^T}{z- \lambda_i} \right)} \times ||E|| \times \Norm{\left(\sum_{i >r} \frac{u_iu_i^T}{z- \lambda_i} \right)} dz\\
& \leq \int_{\Gamma_1} \frac{1}{\min_{i>r} |z- \lambda_i|} ||E|| \frac{1}{\min_{i>r} |z- \lambda_i|} dz \\
& = ||E|| \int_{\Gamma_1}  \frac{1}{ \min_{i>r} |z-\lambda_i|^2} dz \\
& \leq ||E|| \int_{-T}^{T} \frac{1}{\min_{i>r} ((x_0-\lambda_i)^2+t^2)} dt \\
& \leq ||E|| \int_{-T}^{T} \frac{1}{t^2+ (\lambda_p/4)^2 } dt \\ &(\text{since}\, \norm{x_0 -\lambda_i} \geq \norm{\lambda_p -\delta_p/2 - \lambda_{r+1}} \geq \frac{\norm{\lambda_p}}{2} -\frac{\delta_p}{2} \geq \frac{\norm{\lambda_p}}{4} \,\,\,\forall\, i> r) \\
& \leq \frac{16 ||E||}{\norm{\lambda_p}} \,\,\,(\text{by Lemma \ref{ingegralcomputation1}}).
\end{split}
\end{equation}

\noindent Finally , we consider the last term

\begin{equation} \label{lastterm1}
\begin{split}
 \int_{\Gamma_1} \Norm{ \sum_{\substack{i \leq r < j \\ i > r \geq j}} \frac{1}{(z-\lambda_i)(z-\lambda_j)} u_i u_i^T E u_j u_j^T} dz & \leq 2 \int_{\Gamma_1} \Norm{ \left(\sum_{i \leq r} \frac{u_iu_i^T}{z- \lambda_i} \right) E \left( \sum_{j > r} \frac{u_j u_j^T}{z -\lambda_j} \right)} dz \\
& \leq 2 \int_{\Gamma_1} \Norm{\left(\sum_{i \leq r} \frac{u_iu_i^T}{z- \lambda_i} \right)} \times ||E|| \times \Norm{\left(\sum_{j >r} \frac{u_j u_j^T}{z- \lambda_j} \right)} dz\\
& \leq 2 \int_{\Gamma_1} \frac{1}{\min_{i\leq r} |z- \lambda_i|} ||E|| \frac{1}{\min_{j>r} |z- \lambda_j|} dz \\ 
  & \leq 2 ||E|| \int_{\Gamma_1}  \frac{1}{\min_{i \leq r < j} \norm{(z-\lambda_i)(z-\lambda_j)}} dz \\
 & \leq 2 ||E|| \int_{-T}^{T} \frac{1}{\sqrt{(t^2+(\delta_p/2)^2)(t^2+(\lambda_p/4)^2)}} dt \\
 & \text{since}\,\, |z - \lambda_i| \geq \sqrt{t^2+(\delta_p/2)^2}, |z-\lambda_j| \geq \sqrt{t^2+(\lambda_p/4)^2} \,\, \forall i \leq r < j\\
 & = 4 ||E|| \int_{0}^{T}  \frac{1}{\sqrt{(t^2+(\delta_p/2)^2)(t^2+(\lambda_p/4)^2)}} dt. 
 \end{split}
\end{equation}

\noindent Notice that  by Cauchy-Schwartz inequality, 

\begin{equation} \label{lasterm2}
\begin{split}
\int_{0}^{T}  \frac{dt}{\sqrt{(t^2+(\delta_p/2)^2)(t^2+(\lambda_p/4)^2)}}  & \leq \int_{0}^{T} \frac{2}{(t+\delta_p/2)(t+\norm{\lambda_p}/4)} dt \\
& =\frac{2}{\norm{\lambda_p}/4 -\delta_p/2} \int_{0}^{T} \left(\frac{1}{t+\delta_p/2} -\frac{1}{t+ \norm{\lambda_p}/4} \right)dt \\
& = \frac{2}{\norm{\lambda_p}/4 -\delta_p/2} \left[ \log \left( \frac{T+\delta_p/2}{\delta_p/2} \right) - \log \left( \frac{T+|\lambda_p/4|}{|\lambda_p/4|} \right)  \right]\\
& \leq \frac{16}{\norm{\lambda_p}} \times \log \left( \frac{2T+\delta_p}{\delta_p} \right) \,\,(\text{since}\,|\lambda_p|/4 -\delta_p/2 \geq |\lambda_p|/8).
\end{split}
\end{equation}

Together \eqref{lastterm1} and \eqref{lasterm2} imply that the last term is at most 
$$ \frac{64 ||E||}{\norm{\lambda_p}} \times \log \left( \frac{2T+\delta_p}{\delta_p} \right).$$
These estimations imply that 
\begin{equation} 
M_1 \leq \frac{8r^2x}{\delta_p} + \frac{16 ||E||}{\norm{\lambda_p}} + \frac{64||E||}{\norm{\lambda_p}} \times \log \left( \frac{2T+\delta_p}{\delta_p} \right).
\end{equation}
Since $T=2 \sigma_1 \geq 2 |\lambda_p| \geq 8 \delta_p$, thus $\frac{3T}{\delta_p} \geq 24$ and $\log  \left( \frac{3T}{\delta_p} \right) \geq \log 24 > 3$. We further obtain 
\begin{equation}\label{M_1F1bound}
M_1 \leq  \frac{8r^2x}{\delta_p} + \frac{6||E||}{\norm{\lambda_p}} \log \left( \frac{3T}{\delta_p} \right)+ \frac{64||E||}{\norm{\lambda_p}}  \log \left( \frac{2T+\delta_p}{\delta_p} \right) \leq 70 \left( \frac{||E||}{\norm{\lambda_p}} \log \left( \frac{3T}{\delta_p}\right) + \frac{r^2 x}{\delta_p} \right).
\end{equation}

\subsection{Proof of Lemma \ref{lemma: M_1f2}}

Define $x_1:=\max_{1 \leq i,j \leq r_1} \norm{u_i^T E u_j}$. Using triangle inequality, we have 
\begin{equation}
\begin{split}
M_1 & \leq \int_{\Gamma_1} \Norm{ \sum_{i,j \leq r_1} \frac{z}{(z-\lambda_i)(z-\lambda_j)} u_i u_i^T E u_j u_j^T} dz + \int_{\Gamma_1} \Norm{ \sum_{i,j > r_1} \frac{z}{(z-\lambda_i)(z-\lambda_j)} u_i u_i^T E u_j u_j^T} dz \\
& + \int_{\Gamma_1} \Norm{ \sum_{\substack{i \leq r_1 < j \\ i > r_1 \geq j}} \frac{z}{(z-\lambda_i)(z-\lambda_j)} u_i u_i^T E u_j u_j^T} dz.
\end{split}
\end{equation}

\noindent Consider the first term 
\begin{equation}
\begin{split}
\int_{\Gamma_1} \Norm{ \sum_{i,j \leq r_1} \frac{1}{(z-\lambda_i)(z-\lambda_j)} u_i u_i^T E u_j u_j^T} dz & \leq \sum_{i,j \leq r_1} x_1 \int_{\Gamma_1} \frac{\sqrt{a_0^2+t^2}}{\sqrt{((a_0 -\lambda_i)^2+t^2)((a_0 -\lambda_j)^2+t^2)}} dt \\
& \leq r_1^2 x_1 \int_{-T}^{T} \frac{a_0+ |t|}{t^2 +(\delta_k/2)^2} dt \\
& =  r_1^2 x_1 \left( \int_{-T}^{T} \frac{a_0}{t^2+(\delta_k/2)^2} dt + \int_{0}^{T} \frac{2t}{t^2+(\delta_k/2)^2} dt \right) .
\end{split}
\end{equation}

By Lemma \ref{ingegralcomputation1}, we have 
$$\int_{-T}^{T} \frac{a_0}{t^2+(\delta_k/2)^2} dt \leq \frac{8a_0}{\delta_k}.$$

Now we consider the second term

\begin{equation}
\begin{split}
\int_{0}^{T} \frac{2t}{t^2+(\delta_k/2)^2} dt & = \int_{(\delta_k/2)^2}^{T^2+ (\delta_k/2)^2} \frac{1}{u} du \,\,\, ( u= t^2+ (\delta_k/2)^2) \\
& = \log \left( \frac{T^2+(\delta_k/2)^2}{(\delta_k/2)^2}   \right) \\
& = \log \left(\frac{4T^2+\delta_k^2}{\delta_k^2} \right) \leq 2 \log \left( \frac{3T}{\delta_k} \right).
\end{split}
\end{equation}
\noindent Therefore,

\begin{equation} \label{firstermf2M1}
 \int_{\Gamma_1} \Norm{ \sum_{i,j \leq r_1} \frac{z}{(z-\lambda_i)(z-\lambda_j)} u_i u_i^T E u_j u_j^T} dz \leq r_1^2 x_1 \left( \frac{8 a_0}{\delta_k} + 2 \log \left( \frac{3T}{\delta_k} \right)\right).
\end{equation}

\noindent Next, applying the  argument from proof of Lemma \ref{F1f1M1bound} (previous subsection), we bound the second term as follows
\begin{equation}
\begin{split}
\int_{\Gamma_1} \Norm{ \sum_{i,j > r_1} \frac{1}{(z-\lambda_i)(z-\lambda_j)} u_i u_i^T E u_j u_j^T} dz & \leq ||E|| \int_{\Gamma_1}  \frac{|z|}{ \min_{i>r_1} |z-\lambda_i|^2} dz \\
& \leq ||E|| \int_{-T}^{T} \frac{|z|}{ \min_{i>r_1} (a_0-\lambda_i)^2+t^2} dt \\
& \leq ||E|| \int_{-T}^{T} \frac{\sqrt{a_0^2+t^2}}{t^2+ (\lambda_k/4)^2 } dt \\ &(\text{since}\, a_0 -\lambda_i \geq \lambda_k -\delta/2 - \lambda_{r+1} \geq \lambda_k/2 -\delta/2 \geq \lambda_k/4) \\
& \leq ||E||\int_{-T}^{T} \frac{a_0+|t|}{t^2+ (\lambda_k/4)^2 } dt
\end{split}
\end{equation}

\noindent Similarly, 

\begin{equation}
\begin{split}
\int_{-T}^{T} \frac{a_0+|t|}{t^2+ (\lambda_k/4)^2 } dt & \leq \frac{16 a_0}{\lambda_k} + \log \left( \frac{T^2 + (\lambda_p/4)^2}{(\lambda_p/4)^2} \right) \\
& \leq \frac{16 a_0}{\lambda_k} + \log \left(\frac{2T^2}{\delta_k^2}  \right) \\
& \leq 16 + 2\log \left(\frac{2T}{\delta_k}  \right) \,\,\,(\text{since}\,\, a_0 \leq \lambda_k).
\end{split}
\end{equation}

\noindent It follows that 
\begin{equation} \label{secondtermf2M1}
\int_{\Gamma_1} \Norm{ \sum_{i,j > r_1} \frac{z}{(z-\lambda_i)(z-\lambda_j)} u_i u_i^T E u_j u_j^T} dz  \leq ||E|| \left(16 + 2\log \left(\frac{2T}{\delta_k}  \right)   \right).
\end{equation}

\noindent Now we consider the last term 

\begin{equation} \label{lastterm1f2}
\begin{split}
 \int_{\Gamma_1} \Norm{ \sum_{\substack{i \leq r_1 < j \\ i > r_1 \geq j}} \frac{1}{(z-\lambda_i)(z-\lambda_j)} u_i u_i^T E u_j u_j^T} dz & \leq 2 ||E|| \int_{\Gamma_1}  \frac{|z|}{\min_{i\leq r_1 <j} \norm{(z-\lambda_i)(z-\lambda_j)}} dz \\
 & \leq 2 ||E|| \int_{-T}^{T} \frac{|z|}{\sqrt{(t^2+(\delta_k/2)^2)(t^2+(\lambda_k/4)^2)}} dt \\
 & = 4 ||E|| \int_{0}^{T} \frac{\sqrt{a_0^2+t^2}}{\sqrt{(t^2+(\delta_k/2)^2)(t^2+(\lambda_k/4)^2)}} dt \\
 & \leq 4 ||E|| \int_{0}^{T} \frac{a_0+t}{\sqrt{(t^2+(\delta_k/2)^2)(t^2+(\lambda_k/4)^2)}} dt.
  \end{split}
\end{equation}

\noindent  By  Cauchy-Schwartz inequality, we have 
\begin{equation} \label{lastterm2f2}
\begin{split}
\int_{0}^{T} \frac{a_0+t}{\sqrt{(t^2+(\delta_k/2)^2)(t^2+(\lambda_k/4)^2)}} dt &  \leq \int_{0}^T \frac{2(a_0+t)}{(t+\delta_k/2)(t+\lambda_k/4)} dt \\
& = \int_{0}^T \frac{2(a_0-\delta_k/2)}{(t+\delta_k/2)(t+\lambda_k/4)} dt + \int_{0}^T \frac{2(\delta_k/2+t)}{(t+\delta_k/2)(t+\lambda_k/4)} dt \\
& = \int_{0}^T \frac{2(a_0-\delta_k/2)}{(t+\delta_k/2)(t+\lambda_k/4)} dt + \int_{0}^{T} \frac{2}{t+\lambda_k/4} dt \\
& \leq \int_{0}^T \frac{2(a_0-\delta_k/2)}{(t+\delta_k/2)(t+\lambda_k/4)} dt + 2 \log \left(\frac{T+\lambda_k/4}{\lambda_k/4} \right) \\
& \leq \int_{0}^T \frac{2(a_0-\delta_k/2)}{(t+\delta_k/2)(t+\lambda_k/4)} dt + 2 \log \left( \frac{2T}{\delta_k} \right).
\end{split}
\end{equation}

\noindent Similar to \eqref{lasterm2}, since $a_0 \leq \lambda_k$, we also have 
\begin{equation} \label{lasterm2f2}
\int_{0}^T \frac{2(a_0-\delta_k/2)}{(t+\delta_k/2)(t+\lambda_k/4)} dt \leq \frac{ 16(a_0 -\delta_k/2)}{\lambda_k} \log \left( \frac{2T+\delta_k}{\delta_k} \right) \leq 16 \log \left( \frac{2T+\delta_k}{\delta_k} \right).
\end{equation}

The estimates  \eqref{lastterm1f2}, \eqref{lastterm2f2} and \eqref{lasterm2f2}  together imply that the last term is at most 
\begin{equation} \label{Lastermf2}
 72 ||E|| \log \left( \frac{2T+\delta_k}{\delta_k} \right) \leq 72 ||E|| \left( \frac{3T}{\delta_k} \right).
\end{equation}

Combining  \eqref{firstermf2M1}, \eqref{secondtermf2M1} and \eqref{Lastermf2}, we finally obtain that 
\begin{equation} \label{M_1F1f2bound}
\begin{split}
M_1 & \leq r_1^2 x_1 \left( \frac{ 8a_0}{\delta_k} + 2\log \left(\frac{3T}{\delta_k} \right)   \right) +  ||E|| \left( 16 + 2 \log \left( \frac{2T}{\delta_k} \right) \right) + 72 ||E|| \log \left( \frac{3T}{\delta_k} \right) \\
& \leq r_1^2 x_1 \left( \frac{ 8a_0}{\delta_k} + 2\log \left(\frac{6 \sigma_1}{\delta_k} \right)   \right) + 80 ||E|| \log \left(\frac{6 \sigma_1}{\delta_k}  \right) \,\,(\text{since}\,\,  6 \log \left( \frac{6 \sigma_1}{\delta_k}\right) \geq 6 \log 24 > 16 )\\
& \leq r^2 \bar{x} \left( \frac{ 8a_0}{\delta_k} + 2\log \left(\frac{6 \sigma_1}{\delta_k} \right)   \right) + 80 ||E||\log \left(\frac{6 \sigma_1}{\delta_k}  \right).
\end{split}
\end{equation}

\subsection{Proof of Lemma \ref{lemma: N_1,3} and Lemma \ref{lemma: N2,4}}
We bound $N_1, N_2$. The treatment of $N_3,N_4$ is similar and omitted.  First, notice that 

\begin{equation} \label{F_1f2inequality2}
\begin{split}
N_1=& \int_{\Gamma_1}  \frac{|z|}{ \min_i |z-\lambda_i|^2} dz \\
& \leq \int_{-T}^{T} \frac{\sqrt{a_0^2+t^2}}{t^2 + (\delta_k/2)^2} dt \,\, (\text{by construction of $\Gamma$ that }\, |a_0 -\lambda_i| \geq \delta_k/2\,\,\forall 1 \leq  i \leq n) \\
&  \leq \int_{-T}^{T} \frac{\norm{a_0}}{t^2 +(\delta_k/2)^2} dt + \int_{-T}^{T} \frac{|t|}{t^2+(\delta_k/2)^2} dt \\
& = \norm{\frac{2 a_0}{\delta_k}} \int_{-2T/\delta_k}^{2T/\delta_k} \frac{1}{t^2+1}dt + \int_{-2T/\delta_k}^{2T/\delta_k} \frac{|t|}{t^2+1} dt \\
& \leq \norm{\frac{8 a_0}{\delta_k}} + 2 \log \left( \norm{\frac{2T}{\delta_k}}^2 +1 \right)  \,\,(\text{by Lemma \ref{ingegralcomputation1}}) \\
& \leq  \norm{\frac{8 a_0}{\delta_k}} + 4 \log \norm{\frac{3T}{\delta_k}}
\end{split}
\end{equation}

\noindent Next, we have 
\begin{equation} \label{F_1f2inequality3}
\begin{split}
N_2 = & \int_{\Gamma_2}  \frac{|z|}{ \min_i |z-\lambda_i|^2 } dz \\
& = \int_{\Gamma_2}  \frac{\sqrt{x^2+T^2}}{ \min_i ((x-\lambda_i)^2+T^2)} dx \\
& \leq \int_{\Gamma_2} \frac{\sqrt{2} T}{T^2} dx \,\,\,(\text{since}\,\,\, x \leq a_1 \leq T) \\
& = \frac{ \sqrt{2} |a_1 - a_0|}{T}.
\end{split}
\end{equation}

\noindent By similar arguments, we can prove
\begin{equation} \label{F_1f2inequality4}
\begin{split}
& N_3 \leq \frac{4 a_1}{a_1 -\lambda_1} + 4 \log \norm{\frac{3T}{a_1 -\lambda_1}}, \\
& N_4 \leq \frac{ \sqrt{2} |a_1 - a_0|}{T}.
\end{split}
\end{equation}

{\it Acknowledgment.} The research is partially supported by Simon Foundation award SFI-MPS-SFM-00006506 and NSF grant AWD 0010308.


\begin{thebibliography}{1}
\providecommand{\url}[1]{\texttt{#1}}
\expandafter\ifx\csname urlstyle\endcsname\relax
  \providecommand{\doi}[1]{doi: #1}\else
  \providecommand{\doi}{doi: \begingroup \urlstyle{rm}\Url}\fi
\bibitem{AMc} D. Achlioptas and F. McSherry. Fast computation of low rank matrix approximations. \textit{In
Proceedings of the thirty-third annual ACM symposium on Theory of computing}, 611–618. ACM, 2001.
\bibitem{AFKMcS1} Y. Azar, A. Flat, A. Karlin, F. McSherry, and J. Saia, \textit{Spectral analysis of data}, Proceedings of the thirty-third annual ACM symposium on Theory of computing. (2001), 619-626. doi:10.1145/380752.380859.
\bibitem{B-GN1} F. Benaych-Georges and R.R. Nadakuditi, \textit{ The eigenvalues and eigenvectors of finite, low rank perturbations of large random matrices}, Adv. Math. \textbf{227} (2011), no. 1, 494-521.
\bibitem{Book1} R. Bhatia. \textit{Matrix Analysis}, Springer: New York, NY, 2013.

\bibitem{BDDMR1} R. Bhattacharjee, G. Dexter, P. Drineas, C. Musco, and A. Ray, \textit{Sublinear time eigenvalue approximation vua random sampling,} Algorithmica. \textbf{86}, no.6 (2024), 1764-1829. 10.1007/s11042-023-15819-7.


\bibitem{BCN1} C. Bordenave, S. Coste, and R. Nadakuditi, \textit{Detection thresholds in very sparse matrix completion}, Foundations of Computational Mathematics. \textbf{23}, no.5 (2022), 1619-1743. 10.1007/s10208-022-09568-6.




\bibitem{CmFmSq1} J. Chavarria-Molina, J. Fallas-Monge, and P. Soto-Quiros, \textit{Effective implementation to reduce execution time of a low-rank matrix approximation problem,} Journal of Computational and Applied Mathematics. \textbf{401}, no.C (2022). 10.1016/j.cam.2021.113763.

\bibitem{CCF1} Y. Chen, C. Chen, and J. Fan, \textit{ Asymmetry helps: Eigenvalue and Eigenvector analyses of asymmetrically perturbed low-rank matrices}, Ann. Stat. \textbf{49} (2021), no. 1, 435.
\bibitem{CLX1} X-S Chen, W. Li and W.W. Xu, \textit{Perturbation analysis of the eigenvector matrix and singular vector matrices}, Taiwanese Journal of Mathematics. \textbf{16}, No.1, (2012), 179-194.
\bibitem{CLK1} X. Chen, M. Lyu, and I. King, \textit{Toward efficient and accurate covariance matrix estimation on compressed data}, Proceedings of the 34th International Conference on Machine Learning. \textbf{70} (2017), 767-776. 10.5555/3305381.3305381. 
\bibitem{CYZLK1} X. Chen, H. Yang, S. Zhao, M. Lyu, and I. King, \textit{Effective Data-Aware covariance estimator from compressed data,} IEEE Transactions on Neutral Networks and Learning Systems. (2020), 1-14/ 10.1109/TNNLS.2019.2929106.


\bibitem{NAbook} E. Darve and M. Wootters, \textit{Numerical Linear Algebra with Julia,} SIAM Philadelphia, 2021. 

\bibitem{DKoriginal} C. Davis, W.M. Kahan. The rotation of eigenvectors by a perturbation. \textit{III. SIAM J. Numer. Anal.}7, 1970, pp 1-46.

\bibitem{DOT1} J. Draisma, G. Ottaviani, A. Tocino, \textit{
Best rank-k approximations for tensors: generalizing Eckart-Young,} Res. Math. Sci. \textbf{5} (2018), no.2, Paper No. 27, 13 pp.
\bibitem{DriM1} P. Drines and M.W. Mahoney, \textit{Randomization offers new benefits for large-scale linear algebra computations}, RandNLA: randomizerd numerical linear algebra. doi:10.1145/2842602.

\bibitem{DM1} F.J. Dyson, M. Lal Mehta. Random matrices and the statistical theory of energy levels IV. \textit{J. math. Phys}, 4:701–12, 1963.


\bibitem{EBW1} J. Eldridge, M. Belkin and Y. Wang, \textit{Unperturbed: spectral analysis beyond Davis-Kahan}, Algorithmic learning theory. PMLR (2018), 321-358.
\bibitem{G1} J. Ginibre. Statistical ensembles of complex, quaternion, and real matrices. \textit{Journal of Mathematical Physics}, 6(3):440–449, 1965.
\bibitem{G2} V.L. Girko. \textit{Circular law. Theory of Probability \& Its Applications}, 29(4):694–706, 1985.
\bibitem{GSS} F. Gotze, H. Sambale, A. Sinulis. Concentration inequalities for polynomials in $\alpha$-sub-exponential random variables. ArXiv:1903.05964v1.
\bibitem{GTV1} L. Grubišić, N. Truhar, and S. Miodragović, \textit{
The rotation of eigenspaces of perturbed matrix pairs II}, Linear Multilinear Algebra. \textbf{62}(2014), no.8, 1010–1031.
\bibitem{HLMNV1} W. Hachem, P. Loubaton, X. Mestre, J. Najim, P. Vallet, \textit{A subspace estimator for fixed rank perturbations of large random matrices}, J. Multivariate Anal. \textbf{114} (2013), 427–447.


\bibitem{H1} W. Hoeffding. Probability inequalities for sums of bounded random variable. \textit{Journal of the American Statistical Association}, 58 issue 301, 13-30, 1963.
\bibitem{HJBook} R.A. Horn, and C.R. Johnson, \textit{Matrix Analysi}, Cambridge University Press, 2012.

\bibitem{JW1} M. Jirak and M. Wahl, \textit{Perturbation bounds for eigenspaces under a relative gap condition}, Proc. Amer. Math. Soc. \textbf{148} (2020), no.2, 479–494.

\bibitem{Kato1} T. Kato, \textit{Perturbation Theory for Linear Operators}, Classics in Mathematics, Springer: New York, NY, 1980.
\bibitem{KX1} V. Koltchinskii and D. Xia, \textit{ Perturbation of linear forms of singular vectors under Gaussian noise}, High dimensional probability VII, Progr. Probab. \textbf{71} (2016), Springer, 397-423.
\bibitem{KL1} V. Koltchinskii, K. Lounici, \textit{Concentration inequalities and moment bounds for sample covariance
operators}, Bernoulli. \textbf{23}, (2017), 110–133.
\bibitem{Li1} R-C. Li, \textit{Relative Perturbation Theory: (II) Eigenspace and Singular Subspace Variations}, SIAM Journal on Matrix Analysis and Applications. \textbf{20}(2), 471-492, 1998.
\bibitem{Vishnoi1} O. Mangoubi, N.K. Vishnoi. Re-analyze Gauss: Bounds for private matrix approximation via Dyson Brownian motion. \textit{NeurIPS}, 2022.
\bibitem{Vishnoi2} O. Mangoubi, N.K. Vishnoi. Private Covariance Approximation and Eigenvalue-Gap Bounds for Complex Gaussian Perturbations. ArXiv:2306.16648v1, 2023.

\bibitem{MSZ1} S. Morozov, M. Smirnov, N. Zamarashkin, \textit{
On the optimal rank-1 approximation of matrices in the Chebyshev norm}, Linear Algebra Appl. \textbf{679} (2023), 4–29.
\bibitem{PK1} D. Persson and D. Kressner, \textit{
Randomized low-rank approximation of monotone matrix functions},
SIAM J. Matrix Anal. Appl. \textbf{44} (2023), no.2, 894–918.

\bibitem{OVK 13} S. O'Rourke, V. Vu, and K. Wang, \textit{Random Perturbation of low rank matrices: improving classical bounds}, Linear Algebra and its Applications. \textbf{540}, 2018, 26-59. 

\bibitem{OVK 22} S. O'Rourke, V. Vu, and K. Wang, \textit{Matrices with Gaussian noise: Optimal estimates for singular subspace perturbation}, IEEE Transactions on Information Theory. (2023).

\bibitem{RV1} M. Rudelson, R. Vershynin, \textit{Sampling from large matrices: an approach through geometric functional analysis}, J. ACM. \textbf{21}, (2007), 19 pp.
\bibitem{RV2} M. Rudelson, R. Vershynin, \textit{Small ball probabilities for linear images of high-dimensional distributions,} International Mathematics Research Notices 2015. \textbf{19}, (2015), 9594-9617.
\bibitem{SKLZ1} A. Saade, F. Krzakala, M. Lelarge, and L. Zdeborova, \textit{Fast Randomized Semi-Supervised clustering}, Journal of Physics: Conference Series. \textbf{1036} (2018). 10.1088/1742-6596/1036/1/012015. 
\bibitem{SS} W. Schudy, M. Sviridenko. Concentration and Moment Inequalities for Polynomials of Independence Random Variables.  ArXiv:1104.4997v3.
\bibitem{SN1} T.A.B. Snijders, K. Nowicki, Estimation and prediction for stochastic blockmodels for graphs with latent block structure. \textit{Journal of Classification}, \textbf{14}(1), 1997, pp. 75-100.
\bibitem{SWZC1} Z. Song, D. Woodruff, P. Zhong, and T. Chan, \textit{Relative error tensor low rank approximation,} Proceedings of the 30th Annual ACM-SIAM Symposium on Discrete Algorithms. (2019), 2772-2789. 10.5555/3310435.3310607. 

\bibitem{SqCFT1} P. Soto-Quiros, J. Chavarria-Molina, J. Fallas-Monge, and A. Torokhti, \textit{Fast multiple rank-constrained matrix approximation}, SeMA Journal. (2023) 10.1007/s40324-023-00340-6.

\bibitem{St1} W.L. Steiger. Some Kolmogoroff-type inequalities for bounded random variables. \textit{Biometrika}, 54, 641-647, 1967.


\bibitem{SS1} G.W. Stewart, J. Sun. \textit{Matrix Perturbation Theory}. Academic Press. 
\bibitem{Sun1} J-G. Sun, \textit{Perturbation bounds for eigenspaces of a definite matrix pair}, Acta Math. Sinica. \textbf{24} (1981), no.6, 892-903.

\bibitem{TVV} P. Tran, N. Vishnoi, V. Vu, paper in preparation. 



\bibitem{Vbook} R. Vershynin, \textit{High-Dimensional probability - An introduction with applications in data science}, Cambridge Series in Statistical and Probabilistic Mathematics. Cambridge University Press, New York, 2018. DOI: 10.1017/9781108231596. 

\bibitem{Vu0} V. Vu, \textit{Spectral norm of random matrices}, Combinatorica. \textbf{27}, no.6 (2007), 721-736.
\bibitem{Vu1} V. Vu, \textit{Singular vectors under random perturbation}, Random Structures and Algorithm. \textbf{39} (2011), no. 4, 526-538.


\bibitem{WMHZ1} Y. Wan, A. Ma, W. He, and Y. Zhong, \textit{Accurate Multiobjective low-rank and sparse model for hyperspectral image denoising method,} IEEE Transactions on Evolutionary Computation. \textbf{27}, no. 1 (2023), 37-51. 10.1109/TEVC.2021.3078478.

\bibitem{Wa1} R. Wang. Singular vector perturbation under Gaussian noise. \textit{SIAM J. Matrix Anal. Appl.}, 36(1): 158-177, 2015.

\bibitem{Wei1} M. Wei, \textit{Perturbation theory for the Eckart-Young-Mirsky theorem and the constrained total least squares problem,} Linear Algebra Appl. \textbf{280} (1998), no.2-3, 267–287.

\bibitem{We1} Hermann. Weyl, \textit{Das asymptotische Verteilungsgesetz der Eigenwerte linearer partieller Differentialgleichungen}, Mathematische Annalen. \textbf{71} (1912), no. 4, 441-479.
\bibitem{Wu1} Q. Wu, \textit{Relative perturbation bounds for the eigenspace of generalized extended matrices under multiplicative perturbation}, Math. Theory Appl. \textbf{32} (2012), no.2, 53-59.

\bibitem{XY1} D. Xia and M. Yuan, \textit{Effective Tensor Sketching via Sparsification,} IEEE Transactions on Information Theory. \textbf{67}, no. 2 (2021), 1356-1369. 10.1109/TIT.2021.3049174.
\bibitem{YZ1} R. Yuster, U. Zwick, \textit{Fast sparse matrix multiplication}, ACM Transactions on Algorithms (TALG). \textbf{1} (2005), Issue 1, 2 - 13 https://doi.org/10.1145/1077464.1077466.

\bibitem{Z1} Y. Zhong, \textit{Eigenvector Under Random Perturbation: A Nonasymptotic Rayleigh-Schr\"{o}dinger Theory}, (2017). arXiv preprint arXiv:1702.00139.



 
\end{thebibliography}
\end{document}